\documentclass[a4paper,reqno,10pt]{amsart}

\usepackage[USenglish]{babel}

\usepackage{bbm}
\def\1{\mathbbm{1}}
\usepackage{amsmath,amssymb}
\usepackage{amsmath,amsthm}
\usepackage[mathscr]{euscript}
\usepackage{mathrsfs}
\usepackage{chngcntr}
\usepackage{amsfonts}
\usepackage{amssymb}
\usepackage{latexsym}
\usepackage{tikz}
\usepackage{graphics}

\usepackage{verbatim}
\usepackage{stackrel}

\usepackage[shortlabels]{enumitem}
	
\usepackage{comment,cite}	
\usepackage{hyperref}
\hypersetup{
    colorlinks=true,
    linkcolor=blue,
    citecolor=red,
    filecolor=magenta,      
    urlcolor=cyan,
    bookmarks=true,
    linktocpage=true,
}

\usepackage[utf8]{inputenc}

\usepackage[shortlabels]{enumitem}
\usepackage{marginnote}
\usepackage{color}

\usepackage{stmaryrd}
\usepackage{tikz}
\usepackage{tikz-cd}

\newtheorem{definition}{Definition}[section]
\newtheorem{thm}[definition]{Theorem}
\newtheorem{proposition}[definition]{Proposition}
\newtheorem{lemma}[definition]{Lemma}
\newtheorem{corollary}[definition]{Corollary}
\newtheorem{example}[definition]{Example}
\newtheorem{remark}[definition]{Remark}

\numberwithin{equation}{section}

\begin{document}

\title[On the admissibility of observation for evolution families]{On the admissibility of observation operators for evolution families}
\author{Yassine Kharou}
\address{Yassine Kharou, Department of Mathematics, Faculty of Sciences, Ibn Zohr University, Hay Dakhla, B.P. 8106, 80000 Agadir, Morocco}
\email{yassine.kharou@edu.uiz.ac.ma}


\date{September 20, 2021}

\begin{abstract}
This paper is concerned with unbounded observation operators for non-autonomous evolution equations. Fix $\tau > 0$ and let $\left(A(t)\right)_{t \in [0,\tau]} \subset \mathcal{L}(D,X)$, where $D$ and $X$ are two Banach spaces such that $D$ is continuously and densely embedded into $X$. We assume that the operator $A(t)$ has maximal regularity for all $t \in [0,\tau]$ and that $A(\cdot) : [0,\tau] \to \mathcal{L}(D,X)$ satisfies a regularity condition (viz. relative $p$-Dini for some $p \in (1,\infty)$). At first sight, we show that there exists an evolution family on $X$ associated to the problem $$ \dot{u}(t) + A(t) u(t) = 0 \quad t\text{ a.e. on } [0,\tau], \qquad u(0) = x \in X. $$ Then we prove that an observation operator is admissible for $A(\cdot)$ if and only if it is admissible for each $A(t)$ for all $t \in [0,\tau)$.
\end{abstract}

\maketitle

\section{Introduction}\label{sec:1}

Throughout this paper, $(X,\|\cdot\|)$, $(D,\|\cdot\|_D)$ and $(Y,\|\cdot\|_Y)$ are three Banach spaces such that $D$ is continuously and densely embedded into $X$.

We are interested in the concept of admissible observation operators for non-autonomous linear systems, see, e.g., \cite{HHO19,H06,S02}. Before going into details, let us first recall the definition of such operators in the autonomous case. Let $A:D\to X$ be the generator of a strongly continuous semigroup $\mathbb{T} := (\mathbb{T}(t))_{t\ge 0}$ on $X$, $C \in \mathcal{L}(D,Y)$ and $\theta \in (1,\infty)$. The operator $C$ is called a $\theta$-admissible observation operator for $A$ (or $(C,A)$ is $\theta$-admissible) if for some (hence all) $\alpha>0,$ there exists a constant $\gamma>0$ such that
\begin{align}\label{CT-estim}
\int^{\alpha}_{0} \|C \mathbb{T}(t)x\|^{\theta}_{Y} dt \leq \gamma^{\theta} \|x\|^{\theta} \qquad (x \in D).
\end{align}
More details on admissible observation operators can be found in \cite{JP04}, \cite{L03}, \cite{W89} and \cite[Chapter 4]{TW09}. 

Fix $\tau > 0$ and consider a family of unbounded operators $\left(A(t)\right)_{t \in [0,\tau]} \subset \mathcal{L}(D,X)$ such that $A(\cdot):[0,\tau]\to\mathcal{L}(D,X)$ is strongly measurable and bounded and each $A(t)$ has maximal regularity for all $t \in [0,\tau]$. Now, consider the evolution system
\begin{align}\label{nCP}
\begin{cases}
\dot{u}(t)+A(t)u(t)=0 \qquad  t\text{ a.e. on } [s,\tau],\cr u(s)=x ,
\end{cases}
\end{align}  
where $s \in [0,\tau)$ and $x \in X$. We assume that for each $t \in [0,\tau]$, the operator $A(t)$ has maximal regularity and $A(\cdot)$ satisfies the relative $p$-Dini condition (see section \ref{sec:2} and section \ref{sec:3} for definitions). Then we show in Theorem \ref{Kharou-thm1} that there exists an evolution family $U = (U(t,s))_{0 \leq s \leq t \leq \tau}$ on $X$ which solves \eqref{nCP} and $U(t,s)x \in D$ for all $x \in X$ and for almost every $t \in [s,\tau]$. Under the conditions of Theorem \ref{Kharou-thm1}, an operator $C \in \mathcal{L}(D,Y)$ is a $\theta$-admissible observation operator for $A(\cdot)$ (or $(C,A(\cdot))$ is $\theta$-admissible) if and only if for any $s\in [0,\tau)$, there exists a constant $\gamma>0$ such that
\begin{equation}\label{admissibleobservation}
\int_{s}^{\tau} \| C U(t,s)x \|_{Y}^{\theta} \, dt \leq \gamma^{\theta} \| x \|^{\theta} \qquad (x \in D).
\end{equation}
One may ask if there is a relation between the $\theta$-admissibility of $C$ for $A(\cdot)$ and the $\theta$-admissibility of $C$ for each $A(t)$ for any $t \in [0,\tau]$. If $A(\cdot)$ is H\"older continuous of order $\alpha \in (0,1)$ and each $A(t)$ has maximal regularity, then we show in Theorem \ref{Admissibility-Holder} that $\left(C,A(t)\right)$ is $\theta$-admissible for all $t \in [0,\tau]$ implies that $\left( C, A\left(\cdot\right) \right)$ is $\theta$-admissible for every $\theta \in (1,\frac{1}{1-\alpha})$. By the same arguments, we obtain the implication for every $\theta \in (1,\infty)$ assuming that $A(\cdot)$ is Lipschitz continuous. Furthermore, we prove in Theorem \ref{Kharou-thm2} that $\left( C, A\left(\cdot\right) \right)$ is admissible, if and only if $\left(C,A(t)\right)$ is admissible for any $t \in [0,\tau)$, assuming that $A(\cdot): [0,\tau] \to \mathcal{L}(D,X)$ satisfies the condition {\bf(H)} (see section \ref{sec:3}). 

This paper is organized as follows: In the next section, we provide some prerequisites about the concept of maximal regularity. In Section \ref{sec:3}, we prove the main results of the paper.

\section{Maximal regularity and evolution families}\label{sec:2}
In the following definition, we recall what the property of maximal regularity means for a single operator $A \in \mathcal{L}(D,X)$.
\begin{definition}
We say that $A \in \mathcal{L}(D,X)$ has $L^p$-maximal regularity ($p \in (1,\infty)$) and we write $A \in \mathscr{MR}_p$ if for all bounded interval $(a,b) \subset \mathbb{R}$ $(a<b)$ and all $f \in L^{p}(a,b;X)$, there exists a unique $u \in W^{1,p}(a,b;X) \cap L^{p}(a,b;D)$ such that 
\begin{equation}\label{inhomogeneous cp}
\dot{u}(t) + A u(t) = f(t) \quad t\text{ a.e. on } [a,b], \qquad u(a) = 0.
\end{equation}
\end{definition}
\noindent If $A \in \mathscr{MR}_{p}$ $(p \in (1,\infty))$ then $- A$, seen as an unbounded operator on $X$, generates an analytic semigroup on $X$ \cite{D00,KW04}. The converse is true if $X$ is a Hilbert space \cite[Corollary 1.7]{KW04}. According to \cite{D00},  the property of maximal $L^p$-regularity is independent of the bounded interval $(a,b)$, and if $A \in \mathscr{MR}_{p}$ for some $p \in (1,\infty)$ then $A \in \mathscr{MR}_{q}$ for all $q \in (1,\infty)$. Hence, we can write $A \in \mathscr{MR}$ for short. 

Let $p \in (1,\infty)$ and define the maximal regularity space
\begin{align*}
{\rm MR}_{p}(a,b) := W^{1,p}(a,b;X) \cap L^{p}(a,b;D).\end{align*}
It is equipped with the following norm
\begin{align*} \| u \|_{\mathrm{MR}_{p}(a,b)} :=  \| u \|_{W^{1,p}(a,b;X)} + \| u \|_{L^{p}(a,b;D)} \qquad \left( u \in {\rm MR}_{p}(a,b) \right),\end{align*} or, with the equivalent norm $$ \left( \| u \|_{W^{1,p}(a,b;X)}^p + \| u \|_{L^{p}(a,b;D)}^p \right)^{1/p}.$$ The space ${\rm MR}_{p}(a,b)$ is a Banach space when equipped with one of the above norms. Moreover, we consider the trace space defined by $$ \mathrm{Tr}_{p} := \left \{ u(a) : u \in {\rm MR}_{p}(a,b) \right \},$$ and endowed with the norm $$\|x\|_{\mathrm{Tr}_{p}} := \inf \left\{ \| u \|_{{\rm MR}_{p}(a,b)} : x = u(a) \right\}.$$ The space $\mathrm{Tr}_{p}$ is isomorphic to the real interpolation space $\left(X,D\right)_{1-1/p,p}$ \cite[Chapter 1]{L95}. In particular, $\mathrm{Tr}_{p}$ does not depend on the choice of the interval $(a,b)$ and $ D \hookrightarrow_{d} \mathrm{Tr}_{p} \hookrightarrow_{d} X $. We also note that $$ {\rm MR}_{p}(a,b) \hookrightarrow_{d} C \left( [a,b];\mathrm{Tr}_{p} \right) ,$$ and the constant of the embedding does not depend on the choice of the interval $(a,b)$ \cite[Chapter 3]{A95}.\\
It is often practical and convenient to formulate the property of maximal regularity in terms of the invertibility of the sum of two operators (see, e.g., \cite{D00}). This method is commonly known as the operator sum method. For this, consider the operators $\mathcal{A}$ and  $\mathcal{B}$ on $L^{p}(a,b;X)$ given by: 
\begin{align*}
 & D ( \mathcal{B} ) := \left\{ u \in W^{1,p}(a,b;X) : u(a) = 0 \right\},  \quad \mathcal{B}u := \dot{u} \\
 & D ( \mathcal{A} ) := L^{p}(a,b;D),  \quad (\mathcal{A}u)(s) := A u(s) \; \left(s \in (a,b)\right).
\end{align*}
Now consider the operator
\begin{align*}
 & D ( L_A ) := \left\{ u \in \mathrm{MR}_{p}(a,b) : u(a) = 0 \right\},  \quad L_A := \mathcal{A} + \mathcal{B}.
\end{align*}
It is known that $A \in \mathscr{MR}$ if and only if the operator $L_A$ is invertible and $L_{A}^{-1} \in \mathcal{L}\left(L^{p}(a,b;X),\mathrm{MR}_{p}(a,b) \right)$. Moreover, for $f \in L^{p}(a,b;X)$, the function $u = L_{A}^{-1} f$ is the unique solution in $\mathrm{MR}_{p}(a,b)$ of \eqref{inhomogeneous cp} and $$ \| u \|_{\mathrm{MR}_{p}(a,b)} \leq c \| f \|_{L^{p}(a,b;X)},$$ where $c := \| L_{A}^{-1} \|_{\mathcal{L}\left(L^{p}(a,b;X),\mathrm{MR}_{p}(a,b)\right)} $.\\
We have the following proposition. 
\begin{proposition}\label{maximal regularity intial value in Tr}
Assume that $A \in \mathscr{MR}$. Then for every $f \in L^{p}(a,b;X)$ and every $x \in \mathrm{Tr}_{p}$ there exists a unique $u \in \mathrm{MR}_{p}(a,b)$ satisfying 
\begin{equation}\label{cp initial value in Tr}
\dot{u}(t) + A u(t) = f(t) \quad t\text{ a.e. on } [a,b], \qquad u(a) = x ,
\end{equation}
and $$ \| u \|_{\mathrm{MR}_{p}(a,b)} \leq 2 ( c + 1 ) \left[ \| x \|_{\mathrm{Tr}_{p}} + \| f \|_{L^{p}(a,b;X)} \right]. $$

\end{proposition}

\begin{proof}
\emph{Existence}: Let $f \in L^{p}(a,b;X)$ and $x \in \mathrm{Tr}_{p}$. Then there exists $w \in \mathrm{MR}_{p}(a,b)$ such that $x = w(a)$. By maximal regularity, there exists a unique $v \in \mathrm{MR}_{p}(a,b)$ solution of
\begin{equation}
\dot{u}(t) + A u(t) = - \dot{w}(t) - A w(t) + f(t) \quad t\text{ a.e. on } [a,b], \qquad u(a) = 0.
\end{equation}
Putting $u := v + w$, we have $u \in \mathrm{MR}_{p}(a,b)$ and $u$ satisfies \eqref{cp initial value in Tr}.

\emph{Uniqueness}: Let $u_1,u_2 \in \mathrm{MR}_{p}(a,b)$ be two solutions of \eqref{cp initial value in Tr}. Define $v= u_{1} - u_{2}$. Then $v \in \mathrm{MR}_{p}(a,b)$ is a solution of \eqref{cp initial value in Tr} for the right-hand side $f=0$ and the initial value $x=0$. Therefore, by maximal regularity, $v=0$.

It remains to show the estimate. Let $f \in L^{p}(a,b;X)$ and $x \in \mathrm{Tr}_{p}$. There exists $w \in \mathrm{MR}_{p}(a,b)$ such that $x = w(a)$ and $\| w \|_{\mathrm{MR}_{p}(a,b)} \leq 2 \| x \|_{\mathrm{Tr}_{p}}$. The function $u$ given by $u = v + w$ is the unique solution of the problem \eqref{cp initial value in Tr} and
\begin{align*}
\| u \|_{\mathrm{MR}_{p}(a,b)} &= \| v + w \|_{\mathrm{MR}_{p}(a,b)} \\
 &\leq \| w \|_{\mathrm{MR}_{p}(a,b)} + \| v \|_{\mathrm{MR}_{p}(a,b)} \\
 &\leq 2 \| x \|_{\mathrm{Tr}_p} + c \| - \dot{w} - Aw + f \|_{L^{p}(a,b;X)} \\
 &\leq 2 \| x \|_{\mathrm{Tr}_p} + c \left( \| \dot{w} \|_{L^{p}(a,b;X)} + \| Aw \|_{L^{p}(a,b;X)} \right) + c \| f \|_{L^{p}(a,b;X)} \\
 &\leq 2 \| x \|_{\mathrm{Tr}_p} + c \| w \|_{\mathrm{MR}_{p}(a,b)} + c \| f \|_{L^{p}(a,b;X)} \\
 &\leq 2 ( c + 1 ) \left[ \| x \|_{\mathrm{Tr}_p} + \| f \|_{L^{p}(a,b;X)} \right].
\end{align*}
\end{proof}
The following theorem is well known (see, e.g., \cite{D00}). Here we give another simple proof.
\begin{thm}\label{resolvent estimate}
Assume that $A \in \mathscr{MR}$. Then there exists $M \geq 0$, $r > 0$ such that $\left\{ \lambda \in \mathbb{C} : \mathrm{Re} (\lambda) \geq r \right\} \subset \rho(- A)$, where $\rho(- A)$ is the resolvent set of $-A$,
\begin{equation}
\left\| \left( \lambda + A \right)^{-1} \right\|_{\mathcal{L}(X)}\leq\frac{M}{1+\lvert \lambda \rvert}\qquad(\mathrm{Re}(\lambda) \geq r),
\end{equation} 
and
\begin{equation}
\left\| \left( \lambda + A \right)^{-1} \right\|_{\mathcal{L}(X,D)} \leq 2M  \qquad \left( \mathrm{Re}(\lambda) \geq r \right).
\end{equation}

\end{thm}

\begin{proof}
Let $\lambda \in \mathbb{C}$ such that $\mathrm{Re}(\lambda) > 0$ and let $x \in D$. Consider the function $f$ given by $f(t) := \left( \lambda + A \right)x e^{\lambda t}$  for almost every  $t \in [0,1]$. Since $f \in L^{2}(0,1;X)$, then by maximal regularity, there exists a unique $u \in \mathrm{MR}_{2}(0,1)$ satisfying
\begin{equation}\label{cp init value in D}
\dot{u}(t) + A u(t) = f(t) \quad t\text{ a.e. on } [0,1], \qquad u(0) = x.
\end{equation}
By Proposition \ref{maximal regularity intial value in Tr}, there exists a constant $M \geq 0$ such that 
\begin{equation}\label{cst mr}
\| u \|_{\mathrm{MR}_{2}(0,1)} \leq M \left[ \| x \|_{\mathrm{Tr}_p} + \| f \|_{L^{2}(0,1;X)} \right].
\end{equation}
Consider the function $u$ given by $u(t) := e^{\lambda t} x$. The function $u$ satisfies \eqref{cp init value in D} and $u \in \mathrm{MR}_{2}(0,1)$. Thus $u$ is the unique solution of \eqref{cp init value in D}. \\
Let us now compute $\| f \|_{L^{2}(0,1;X)}$ and $\| u \|_{\mathrm{MR}_{2}(0,1)}$ in order to inject them in \eqref{cst mr}, we have:
\begin{align*}
\| f \|_{L^{2}(0,1;X)} = \left( \int_{0}^{1} \| f(t) \|^{2} dt \right)^{1/2} &= \left( \int_{0}^{1} \| \left( \lambda + A \right) x e^{\lambda t} \|^{2} dt \right)^{1/2} \\
 &= \| \left( \lambda + A \right) x \| \| e^{\lambda \cdot} \|_{L^{2}(0,1)} ,
\end{align*} 
and
\begin{align*}
\| u \|_{\mathrm{MR}_{2}(0,1)} &= \| u \|_{L^{2}(0,1;D)} + \| u \|_{L^{2}(0,1;X)} + \| \dot{u} \|_{L^{2}(0,1;X)} \\
 &= \left( \int_{0}^{1} \| e^{\lambda t} x \|_{D}^{2} dt \right)^{1/2} + \left( \int_{0}^{1} \| e^{\lambda t} x \|^{2} dt \right)^{1/2} + \left( \int_{0}^{1} \| \lambda e^{\lambda t} x \|^{2} dt \right)^{1/2} \\
 &= \| e^{\lambda \cdot } \|_{L^{2}(0,1)} \| x \|_{D} + \| e^{\lambda \cdot } \|_{L^{2}(0,1)} \| x \| +  \lvert \lambda \rvert \| e^{\lambda \cdot } \|_{L^{2}(0,1)} \| x \| .
\end{align*}
It follows that 
\begin{equation}
\| x \|_D + (1 + \lvert \lambda \rvert) \| x \| \leq c_i M \sqrt{\frac{2 \mathrm{Re}(\lambda)}{e^{2 \mathrm{Re}(\lambda)} - 1}} \| x \|_D + M \| \left( \lambda + A \right) x \| ,
\end{equation}
where $c_i$ is the constant of the embedding $D \hookrightarrow_{d} Tr $.\\
There exists $r > 0$ such that $$ \sqrt{\frac{e^{2 \mathrm{Re}(\lambda)} - 1}{2 \mathrm{Re}(\lambda)}} \geq 2 c_i M ,$$ for $\mathrm{Re}(\lambda) \geq r$. Then, we obtain for $\lambda \in \mathbb{C}$ with $\mathrm{Re}(\lambda) \geq r$,
\begin{equation}
\frac{1}{2} \| x \|_D + \left( 1 + \lvert \lambda t \rvert \right) \| x \| \leq M \| \left( \lambda + A \right) x \| .
\end{equation}
Therefore, for $\lambda \in \mathbb{C}$ such that $\mathrm{Re}(\lambda) \geq r$, $\lambda + A$ is injective. Since $A \in \mathscr{MR}$ then $A$ is closed and hence $\left( \lambda + A \right)^{-1}$ is closed, which implies that $\left( \lambda + A \right)^{-1}$ is bounded. Finally, $$ \left\| \left( \lambda + A \right)^{-1} x \right\| \leq \frac{M}{1 + \lvert \lambda \rvert} \| x \| \qquad \left( \mathrm{Re}(\lambda) \geq r, x \in X \right) $$ and $$ \left\| \left( \lambda + A \right)^{-1} x \right\|_D \leq 2M \| x \| \qquad \left( \mathrm{Re}(\lambda) \geq r, x \in X \right). $$ 
\end{proof}

Consider the non-autonomous evolution equation
\begin{equation}\label{homogeneous initial value problem}
\dot{u}(t) + A(t) u(t) = 0 \quad t\text{ a.e. on } [s,\tau], \qquad u(s) = x
\end{equation}
for any $s \in [0,\tau)$ and $x \in X$, where $A(\cdot) : [0,\tau] \to \mathcal{L}(D,X)$ is strongly measurable and bounded such that $A(t) \in \mathscr{MR}$ for all $t \in [0,\tau]$. The following definition recalls what the property of maximal regularity means in the non-autonomous case. For more information, we refer to \cite{ACFP07} and \cite{LE13}.

\begin{definition}
Let $p \in (1,\infty)$ and $A(\cdot) : [0,\tau] \to \mathcal{L}(D,X)$ be a bounded and strongly measurable function. We say that $A(\cdot)$ has $L^p$-maximal regularity on the bounded interval $[0,\tau]$ (and we write $A(\cdot) \in \mathscr{MR}_{p}(0,\tau)$), if and only if for all $[a,b]$ $(a<b)$ a sub-interval of $[0,\tau]$ and all $f \in L^{p}(a,b;X)$, there exists a unique $u \in \mathrm{MR}_{p}(a,b)$ such that
\begin{equation}
\dot{u}(t) + A(t) u(t) = f(t) \quad t\text{ a.e. on } [a,b], \qquad u(a) = 0.
\end{equation}
\end{definition}
\noindent Recently, the authors of \cite{ACFP07} proved that $A(\cdot) \in \mathscr{MR}_{p}(0,\tau)$, assuming that each $A(t)$ has maximal regularity for all $t \in [0,\tau]$ and that $A(\cdot)$ is relatively continuous (in the sense of the following definition).
\begin{definition}\label{Relative continuity}
We say that  $A(\cdot):[0,\tau]\to \mathcal{L}(D,X)$ is  relatively continuous if for all $\varepsilon > 0$ there exist $\delta > 0$ and $\eta \geq 0$ such that for all $s,t \in [0,\tau]$, we have: $$ \| A(t) x - A(s) x \| \leq \varepsilon \| x \|_D + \eta  \| x \| \qquad ( x \in D)$$ whenever $\lvert t-s \rvert \leq \delta$.
\end{definition}
\noindent In the non-autonomous case, the concept of evolution family will play the same role as the notion of semigroup in the autonomous case. However, in the non-autonomous case, we do not have a general theory and results like the Hille-Yosida generation theorem. In order to obtain the existence and the regularity of the evolution family on the state space $X$, one should make additional assumptions on $A(\cdot)$, see \cite{AT87} and \cite{S01}. The following definition explains what we mean by an evolution family.
\begin{definition}\label{evolutions family}
A family $U:=(U(t,s))_{0\le s\le t\le \tau}\subset \mathcal{L}(X)$ is an evolution family on $X$ if
\begin{itemize}
  \item [{\rm (i)}] $U(t,s)=U(t,r)U(r,s),\; U(s,s)=I$, for any $s,t,r \in [0,\tau]$ with $0 \leq s \leq r \leq t \leq \tau$,
  \item [{\rm (ii)}] $\left\{ (t,s) \in [0,\tau]^2 : t \geq s \right\} \ni (t,s)\mapsto U(t,s)$ is strongly continuous.
\end{itemize}

\end{definition}
\noindent Evolution families are very efficient in solving non-autonomous evolution equations. For more details, see \cite{CL99}, \cite{EN00}, \cite{P83}. In fact, If $A(\cdot) : [0,\tau] \to \mathcal{L}(D,X)$ is a relatively continuous function such that $A(t) \in \mathscr{MR}$ for any $t \in [0,\tau]$. Thus by \cite[Theorem 2.7]{ACFP07}, we have $A(\cdot) \in \mathscr{MR}_{p}(0,\tau)$ for any $p \in (1,\infty)$. Using \cite[Lemma 2.2]{ACFP07}, then for every $x \in \mathrm{Tr}_p$ and every $(t,s) \in \Delta_{0,\tau} := \left\{ (t,s) \in [0,\tau]^2 : t \geq s \right\}$, we can define $$ U(t,s) x := u(t) ,$$ where $u$ is the unique solution of the problem $$ \dot{u}(t) + A(t) u(t) = 0 \quad t\text{ a.e. on } [s,\tau], \qquad u(s) = x .$$ Furthermore, the family $(U(t,s))_{0\le s\le t\le \tau}$ is an evolution family on $\mathrm{Tr}_p$ \cite[Proposition 2.3]{ACFP07}. If, in addition, each $A(t)$ is accretive for any $t \in [0,\tau]$, the authors of \cite{ACFP07} showed that $U$ extends to a contractive evolution family on $X$. Under the condition {\bf(H)} (see section \ref{sec:3}) on $A(\cdot)$, we show in Theorem \ref{Kharou-thm1} that the evolution family $U$ on $\mathrm{Tr}_{p}$ extends to an evolution family on $X$. For its proof, we need the following lemma.

\begin{lemma}\label{resolvent estimate t}
Let $A(\cdot) : [0,\tau] \longrightarrow \mathcal{L}(D,X)$ be strongly measurable and relatively continuous. Assume that $A(t) \in \mathscr{MR}$ for all $t \in [0,\tau]$. Then there exist $M_0 \geq 0$, $r_0 > 0$, independent of $t \in [0,\tau]$, such that $\left\{ \lambda \in \mathbb{C} : \mathrm{Re} (\lambda) \geq r_0 \right\} \subset \rho(- A(t))$ and
\begin{equation}
\left\| \left( \lambda + A(t) \right)^{-1} \right\|_{\mathcal{L}(X)} \leq \frac{M_0}{1 + \lvert \lambda \rvert} \qquad \left( \mathrm{Re}(\lambda) \geq r_0 \right)
\end{equation} 
for all $t \in [0,\tau]$.

\end{lemma}

\begin{proof}
By the assumption of relative continuity on $A(\cdot)$, we have for every $t \in [0,\tau]$ there exist $\delta_t > 0$ and $\eta_t \geq 0$ such that for every $s \in [t-\delta_t,t+\delta_t]$ and every $x \in D$, $$ \| A(t)x - A(s)x \| \leq \frac{1}{4 M_t} \| x \|_D + \eta_t \| x \|,$$ where $M_t$ is the constant in Theorem \ref{resolvent estimate} associated with $A(t)$.\\
By \cite[Lemma 2.8]{ACFP07}, there exist a partition $0 = \tau_0 < \tau_1 < \cdots < \tau_n = \tau$ and $t_i \in [0,\tau]$, $i \in \left\{0,1,\cdots,n\right\}$, such that $t_i \in [\tau_i,\tau_{i+1}] \subset [t_i - \delta_{t_i},t_i + \delta_{t_i}]$ for all $i \in \left\{0,1,\cdots,n-1\right\}$.\\
Let $x \in X$, $t \in [0,\tau]$ and let $i \in \left\{0,1,\cdots,n-1\right\}$ be such that $t \in [\tau_i,\tau_{i+1}]$. Then for $\lambda \in \mathbb{C}$ such that $\mathrm{Re}(\lambda) \geq  \sup_{0 \leq i \leq n} r_{t_i}$, where $r_{t_i}$ is the constant given in Theorem \ref{resolvent estimate} associated with $A(t_i)$, we have:
\begin{align*}
\left\| \left( A(t) - A(t_i) \right)\left(\lambda + A(t_i)\right)^{-1} x \right\| &\leq \frac{1}{4 M_{t_i}} \left\| \left(\lambda + A(t_i)\right)^{-1} x \right\|_D + \eta_{t_i} \left\| \left(\lambda + A(t_i)\right)^{-1} x \right\| \\
 &\leq \frac{1}{2} \left\| x \right\| + \frac{\eta_{t_i} M_{t_i}}{1 + \lvert \lambda \rvert} \left\| x \right\| .
\end{align*}
We set $r_0 := \max\left\{4 \sup_{0 \leq i \leq n} \left( \eta_{t_i} M_{t_i} \right) ; \sup_{0 \leq i \leq n} r_{t_i}\right\}$. Thus for $\lambda \in \mathbb{C}$ such that $\mathrm{Re}(\lambda) \geq r_0$, we have: $$\left\| \left( A(t) - A(t_i) \right)\left(\lambda + A(t_i)\right)^{-1}  \right\|_{\mathcal{L}(X)} \leq \frac{3}{4} < 1$$ and 
\begin{align*}
\left(\lambda + A(t)\right)^{-1} &= \left( \lambda + A(t_i) + A(t) - A(t_i) \right)^{-1} \\
 &= \left( \lambda + A(t_i) \right)^{-1} \left( I + \left(A(t)-A(t_i)\right)\left(\lambda+A(t_i)\right)^{-1}\right)^{-1}.
\end{align*}
Therefore,
\begin{align*}
\left\| \left(\lambda + A(t)\right)^{-1} \right\|_{\mathcal{L}(X)} &\leq \left\| \left( \lambda + A(t_i) \right)^{-1} \right\|_{\mathcal{L}(X)} \left\| \left( I + \left(A(t)-A(t_i)\right)\left(\lambda+A(t_i)\right)^{-1}\right)^{-1} \right\|_{\mathcal{L}(X)} \\
 &\leq \frac{M_{t_i}}{1+\lvert \lambda \rvert} \sum_{k=0}^{\infty} \left\| \left(A(t)-A(t_i)\right)\left(\lambda+A(t_i)\right)^{-1}\right\|_{\mathcal{L}(X)}^k \\
 &\leq \frac{4 \sup_{0 \leq i \leq n} M_{t_i}}{1+\lvert \lambda \rvert} .
\end{align*}
\end{proof}

\section{Main results}\label{sec:3}

We need the following regularity assumption on the map $A(\cdot)$.
\begin{definition}
Let $p \in (1,\infty)$. The function $A(\cdot):[0,\tau]\to\mathcal{L}(D,X)$ satisfies the relative $p$-Dini condition if there exist $\eta \geq 0$, $\omega : [0,\tau] \to [0,\infty)$ a continuous function with $\omega(0)=0$ and  
 \begin{align}\label{omega-condition} 
 \int_{0}^{\tau} \left( \frac{\omega(t)}{t} \right)^{p} \mathrm{d}t < \infty 
 \end{align} such that for all $x \in D$, $s,t \in [0,\tau]$, we have: $$\| A(t) x - A(s) x \| \leq \omega \left( \lvert t-s \rvert \right) \| x \|_D + \eta  \| x \|.$$
\end{definition}
\noindent We emphasize that if $A(\cdot)$ satisfies the relative $p$-Dini condition, then $A(\cdot)$ is relatively continuous. Remark that if $A(\cdot)$ is H\"older continuous, that is,
\begin{align*}
\|A(t)x-A(s)x\|\leq \lvert t-s \rvert^{\alpha} \|x\|_D \qquad (t,s\in [0,\tau], \, x \in D)
\end{align*}
for some $\alpha \in (0,1)$, then it satisfies the relative $p$-Dini condition for all $p\in (1,\frac{1}{1-\alpha})$.\\
In the following, we need the condition {\bf(H)} given by:
\begin{itemize}
  \item [{\bf(H)}] For any $t\in [0,\tau],$ $A(t)\in\mathscr{MR}$, and $A(\cdot):[0,\tau]\to\mathcal{L}(D,X)$ satisfies the relative $p$-Dini condition.
\end{itemize}
\begin{remark}
Let $A(\cdot): [0,\tau] \to \mathcal{L}(D,X)$ satisfies the condition {\bf(H)}. The function $A(\cdot)$ is strongly measurable and relatively continuous. Then by Lemma \ref{resolvent estimate t} and by the rescaling technique, we assume without loss of generality, that there exists a constant $M_0 \geq 0$ such that
\begin{equation}
\left\| R \left( \lambda,-A(t) \right) \right\|_{\mathcal{L}(X)} = \left\| \left( \lambda + A(t) \right)^{-1} \right\|_{\mathcal{L}(X)} \leq \frac{M_0}{1 + \lvert \lambda \rvert} \qquad \left( \mathrm{Re}(\lambda) \geq 0 , t \in [0,\tau] \right).
\end{equation}
This implies that for every $t \in [0,\tau]$, the operator $- A(t)$ generates an analytic semigroup $\left( e^{- s A(t)} \right)_{s \geq 0}$, satisfying $$ \left\| e^{-s A(t)} \right\| \leq c \qquad (0 \leq s \leq \tau) $$ and $$ \left\| A(t) e^{-s A(t)} \right\| \leq \frac{c}{s} \qquad (0 < s \leq \tau) ,$$ for some constant $c \geq 0$ (see \cite[Theorem 2.5.2]{P83}).
\end{remark}
The following theorem shows that to $A(\cdot)$, we can associate an evolution family on $X$.
\begin{thm}\label{Kharou-thm1}
Assume that $A(\cdot): [0,\tau] \longrightarrow \mathcal{L}(D,X)$ satisfies the condition {\rm\bf(H)} and  let $U:=(U(t,s)_{0\le s\le t\le \tau}$  be the associated evolution family on the trace space $\mathrm{Tr}_p$ for $p \in (1,\infty)$. Then the following assertions hold:
\begin{itemize}
  \item [{\rm (i)}] $U$ extends to a bounded evolution family on $X$.
  \item [{\rm (ii)}]  The function $u$ given by $u(t) := U(t,0)x$ is the unique solution of the problem
\begin{align*}
\dot{u}(t) + A(t) u(t) = 0 \quad t\text{ a.e. on } [0,\tau], \qquad u(0) = x.
\end{align*}
  \item [{\rm (iii)}]  For $q \in (1,\infty),$ and $x \in X$, $v:t\mapsto v(t)= tU(t,0)x \in {\rm MR}_{q}(0,\tau)$ and
  \begin{align}\label{goud-estomate} \| v \|_{{\rm MR}_{q}(0,\tau)} \leq M \| x \|\end{align}
  for a constant $M \geq 0$, depending on $q$ but independent of $x \in X$. Moreover, the function $u$ given by $u(t) := U(t,0)x$ belongs to the space $$ C\left([0,\tau];X\right) \cap L^{q}_{loc}\left((0,\tau];D\right) \cap W^{1,q}_{loc}\left((0,\tau];X\right).$$
\end{itemize}
\end{thm}
\begin{proof}

\begin{enumerate}
\item By \cite[Theorem 2.7]{ACFP07}, $A(\cdot) \in \mathscr{MR}_{q}(0,\tau)$ for every $q \in (1,\infty)$. Let $x \in \mathrm{Tr}_p$ and $s \in [0,\tau)$. Consider the function $f : [s,\tau] \ni t \mapsto \left( A(s) - A(t) \right)e^{-(t-s)A(s)} x$. The inhomogeneous Cauchy problem
\begin{equation}\label{eq s f}
\dot{u}(t) + A(t) u(t) = f(t) \qquad t\text{ a.e. on } [s,\tau], \qquad u(s) = 0
\end{equation}
has a unique solution $u \in \mathrm{MR}_{p}(s,\tau)$ given by: $ u(t) = U(t, s)x - e^{-(t -s)A(s)} x$. In fact, we have:
\begin{align*}
\int_{s}^{\tau} \| f(t) \|^{p} \mathrm{d}t &= \int_{s}^{\tau} \| \left( A(s) - A(t) \right) e^{-(t-s)A(s)} x \|^{p} \mathrm{d}t \\
 &\leq \int_{s}^{\tau} \left( \omega \left( \lvert t-s \rvert \right) \| e^{-(t-s)A(s)}x \|_{D} + \eta \| e^{-(t-s)A(s)}x \| \right)^{p} \mathrm{d}t \\
 &\leq 2^p \int_{s}^{\tau} \left( \omega ( \lvert t-s  \rvert^{p} \frac{c^{p}}{\lvert t-s \rvert^p}  \| x \|^{p} + \eta^{p} c^{p} \| x \|^{p} \right) \mathrm{d}t \\
 &\leq 2^p c^p \int_{0}^{\tau} \left( \frac{\omega(t)}{t} \right)^p \mathrm{d}t \| x \|^{p} + 2^p \eta^{p} c^p \tau \| x \|^{p} \\
 &\leq 2^p c^p \left[ \int_{0}^{\tau} \left( \frac{\omega(t)}{t} \right)^p \mathrm{d}t + \eta^p \tau \right] \| x \|^{p} ,
\end{align*}
then $f \in L^{p}(s,\tau;X)$. By $L^{p}$-maximal regularity, there exists a unique solution $u \in \mathrm{MR}_{p}(s,\tau)$ of (\ref{eq s f}).\\
Since $x \in \mathrm{Tr}_p$, it is easy to see that $$U(\cdot,s)x - e^{-(\cdot -s)A(s)}x \in \mathrm{MR}_{p}(s,\tau)$$ and
\begin{align*}
\frac{\mathrm{d}}{\mathrm{d}t} \left [ U(t,s)x - e^{-(t-s)A(s)}x \right ] &=  - A(t) U(t,s)x + A(s) e^{-(t-s)A(s)}x \\
 &= - A(t) \left [ U(t,s)x - e^{-(t-s)A(s)}x \right ] \\
 &  \qquad \qquad \qquad + \left(- A(t) + A(s)\right) e^{-(t-s)A(s)}x \\
 &= - A(t) \left [ U(t,s)x - e^{-(t-s)A(s)}x \right ] + f(t),
\end{align*} 
Therefore $U(\cdot,s)x - e^{(\cdot -s)A(s)}x$ is the unique solution of (\ref{eq s f}).

\item By $L^{p}$-maximal regularity, it follows that there exists a constant $\kappa > 0$ such that for every $x \in \mathrm{Tr}_p$,
$$ \| U(\cdot,s)x - e^{-(\cdot - s)A(s)}x \|_{\mathrm{MR}_{p}(s,\tau)} \leq \kappa \| f \|_{L^{p}(s,\tau;X)} $$ then
\begin{align}\label{key inequality}
\| U(\cdot,s)x - e^{-(\cdot - s)A(s)}x \|_{\mathrm{MR}_{p}(s,\tau)} \leq 2 c \kappa \left[ \int_{0}^{\tau} \left( \frac{\omega(t)}{t} \right)^p \mathrm{d}t + \eta^p \tau \right]^{1/p} \| x \| .
\end{align}
From the continuous embedding $ \mathrm{MR}_{p}(s,\tau) \hookrightarrow_{d} C([s,\tau];X)$, we obtain that there exists a constant $M \geq 0$, independent of $s$ and $t$, such that $$\| U(t,s)x - e^{-(t - s)A(s)}x \| \leq M \| x \| \qquad \text{ for } (t,s) \in \Delta_{0,\tau} \text{ and } x \in \mathrm{Tr}_p .$$ The density of $\mathrm{Tr}_p$ on $X$ implies that the evolution family $\left(U(t,s)\right)_{0\leq s \leq t \leq \tau}$ extends to an evolution family on $X$.

\item The evolution family $U = \left(U(t,s)\right)_{0\leq s \leq t \leq \tau}$ extends to an evolution family on $X$, which we also denote by $U$ and there exists a constant $M \geq 0$ such that $$ \| U(t,s) x \| \leq M \| x \| $$ for every $x \in X$ and every $(t,s) \in \Delta_{0,\tau}$.\\
For every $x \in \mathrm{Tr}_p$, we obtain by $L^q$-maximal regularity of $A(\cdot)$ that the function $v(t) := t U(t,0)x$ is the unique solution of the non-homogeneous problem $$ \dot{v}(t) + A(t) v(t) = U(t,0)x \quad t\text{ a.e. on } [0,\tau], \qquad v(0) = 0 .$$ And that there exists a constant $\kappa > 0$ such that $$ \| v \|_{\mathrm{MR}_{q}(0,\tau)} \leq \kappa \left\| U(\cdot,0)x \right\|_{L^{q}(0,\tau;X)} \leq M' \| x \|.$$ The density of $\mathrm{Tr}_p$ on $X$ implies that this estimate holds for every $x \in X$. In particular, for every $x \in X$ and every $q \in (1,\infty)$, we have: $$ U(\cdot,0) x \in L_{loc}^{q}\left((0,\tau];D\right) \cap W_{loc}^{1,q} \left((0,\tau];X\right).$$ The claim follows from the definition of $U$.
\end{enumerate}
\end{proof}

\begin{corollary}\label{Main Corollary}
Let $A(\cdot): [0,\tau] \to \mathcal{L}(D,X)$ satisfies the condition {\bf(H)}. Let $q \in (1,\infty)$. Then for every $x \in X$ and every $f \in L^{q}(0,\tau;X)$ the function $u$ given by 
$$
u(t) := U(t,0)x + \int_{0}^{t} U(t,s) f(s) \, \mathrm{d}s
$$ 
is in the space $$ C\left([0,\tau];X\right) \cap L^{q}_{loc}\left((0,\tau];D\right) \cap W^{1,q}_{loc}\left((0,\tau];X\right)$$ and is the unique solution of the problem 
$$
\dot{u}(t) + A(t) u(t) = f(t) \quad t\text{ a.e. on } [0,\tau], \qquad u(0) = x.
$$

\end{corollary}

\begin{example}
Let $q \in (1,\infty)$ and $\Omega$ be an open set in $\mathbb{R}^n$ such that $\partial \Omega$ is bounded and of class $C^2$. Assume that $a_{ij} \in C ( [0,\tau] \times \bar{\Omega} )$ for $i,j=1,\cdots,n$ is bounded and uniformly elliptic, i.e., $$ \sum_{i,j=1}^{n} a_{ij}(t,x) \xi_{i} \xi_{j} \geq \beta \lvert \xi \rvert^2 $$ for some $\beta > 0$ and all $\xi \in \mathbb{R}^n$, $x \in \bar{\Omega}$, $t \in [0,\tau]$. Assume in addition that $$ \max_{i,j=1,\cdots,n} \, \sup_{x \in \Omega} \, \lvert a_{ij}(t,x) - a_{ij}(s,x) \rvert \leq \omega ( t-s ) \qquad ( t,s \in [0,\tau] ) ,$$ where $\omega : [0,\tau] \longrightarrow [0,\infty)$ is a continuous function satisfying $$ \int_{0}^{\tau} \left( \frac{\omega(t)}{t} \right)^{\alpha} \mathrm{d}t < \infty , $$ for some $\alpha \in (1,\infty)$.\\
Define the partial differential operator $\mathcal{A}(t,x,D)$ by $$ \mathcal{A}(t,x,D) u(x) := \sum_{i,j=1}^{n} a_{ij}(t,x) \partial_{i} \partial_{j} u(x) + b_{0}(t,x) u(x) ,$$ where $b_{0} \in L^{\infty} \left( (0,\tau) \times \Omega \right)$ .\\
Let $D:= W^{2,p}(\Omega) \cap W_{0}^{1,p}(\Omega)$ with $1 < p < \infty$ and define for every $t \in [0,\tau]$ the operator $A(t) \in \mathcal{L}\left(D,L^{p}(\Omega)\right)$ by $$ A(t) u := - \sum_{i,j=1}^{n} a_{ij}(t,\cdot) \partial_{i} \partial_{j} u - b_{0}(t,\cdot)u  \qquad ( u \in D ).$$ It follows from \cite[Theorem 8.2]{DHP03} and \cite[Proposition 1.3]{ACFP07} that $A(t) \in \mathscr{MR}$ for every $t \in [0,\tau]$. Moreover, for every $t,s \in [0,\tau]$ and $u \in D$, we have: 
\begin{align*}
\| A(t) u - A(s) u \|_{L^{p}(\Omega)} &= \left\| \sum_{i,j=1}^{n} \left( a_{ij}(t,\cdot) - a_{ij}(s,\cdot) \right) \partial_{i} \partial_{j} u + \left( b_{0}(t,\cdot) - b_{0}(s,\cdot) \right) u  \right\|_{L^{p}(\Omega)} \\
 &\leq \omega \left( t-s \right) \sum_{i,j=1}^{n} \| \partial_{i} \partial_{j} u \|_{L^{p}(\Omega)} + \eta \| u \|_{L^{p}(\Omega)} \\
 &\leq \omega \left( t-s \right) \| u \|_D + \eta \| u \|_{L^{p}(\Omega)} ,
\end{align*}
for some $\eta \geq 0$. Then by Corollary \ref{Main Corollary}, we obtain that for every $u_0 \in L^{p}(\Omega)$ and every $f \in L^{q}(0,\tau;L^{p}(\Omega))$ there exists a unique solution $$ u \in C([0,\tau]; L^{p}(\Omega) ) \cap L^{q}_{loc}((0,\tau];W^{2,p}(\Omega) \cap W_{0}^{1,p}(\Omega)) \cap W^{1,q}_{loc}((0,\tau];L^{p}(\Omega)) $$ of the problem
$$
\left \{
\begin{array}{r c l}
\partial_t u(t)(x) - \mathcal{A}(t,x,D) u(t)(x) &=& f(t)(x) \qquad (t,x)\text{ a.e on } [0,\tau] \times \Omega,\\ 
u(t)(x) \qquad \qquad &=& 0 \qquad (t,x)\text{ a.e on } [0,\tau] \times \partial \Omega,\\
u(0)(x) \qquad \qquad &=& u_{0}(x) \qquad x\text{ a.e on } \Omega.

\end{array}
\right . 
$$

\end{example}

The following theorem shows a condition on $A(\cdot)$ for which $\left(C,A(t)\right)$ is admissible for all $t \in [0,\tau]$ implies that $\left( C, A\left(\cdot\right) \right)$ is admissible.

\begin{thm}\label{Admissibility-Holder}
Let $A(\cdot): [0,\tau] \to \mathcal{L}(D,X)$ be a H\"older continuous function, i.e.,
\begin{align*}
\|A(t)x-A(s)x\|\leq \lvert t-s \rvert^{\alpha} \|x\|_D \qquad \left( t,s\in [0,\tau], x \in D \right),
\end{align*}
for some $\alpha \in (0,1)$ and let $C \in \mathcal{L}(D,Y)$. Assume that $A(t)\in\mathscr{MR}$ for any $t \in [0,\tau]$. \\
If $\left(C,A(t)\right)$ is $p$-admissible for all $t \in [0,\tau]$, then $\left( C, A\left(\cdot\right) \right)$ is $p$-admissible for every $p < \frac{1}{1-\alpha}$.

\end{thm}

\begin{proof}
Let $0 \leq s < t \leq \tau$, $x \in D$ and $p \in (1,\infty)$. Consider the function $v:[s,t] \ni r \mapsto e^{-(t-r)A(s)} U(r,s)x$. Then $v \in W^{1,p}(s,t;X)$ and
\begin{align*}
\frac{d}{dr} v(r) &= A(s) e^{-(t-r)A(s)} U(r,s) x - e^{-(t-r)A(s)} A(r) U(r,s) x \\
 &= e^{-(t-r)A(s)} ( A(s) - A(r) ) U(r,s) x .
\end{align*}
By integrating between $s$ and $t$, we obtain $$ U(t,s) x = e^{-(t-s) A(s)} x + \int_{s}^{t} e^{-(t-r)A(s)} ( A(s) - A(r) ) U(r,s) x dr .$$ Assume now that $\left(C,A(s)\right)$ is $p$-admissible. Thus
\begin{align}\label{estimate}
\begin{split}
\int_{s}^{\tau} \| C U(t,s) x \|^p dt &\leq 2^p \int_{s}^{\tau} \| C e^{-(t-s) A(s)} x \|^p dt \\
 &\qquad + 2^p \int_{s}^{\tau} \left\| C \int_{s}^{t}  e^{-(t-r)A(s)} ( A(s) - A(r) ) U(r,s) x dr \right\|^p dt \\
 &\leq \gamma_{\tau}^p \|x\|^p + k_{\tau} \int_{s}^{\tau} \left\| ( A(s) - A(r) ) U(r,s) x \right\|^p dr \\
 &\leq \gamma_{\tau}^p \|x\|^p + k_{\tau} \int_{s}^{\tau} \frac{1}{(r-s)^{(1-\alpha)p}} dr \|x\|^p .
 \end{split}
\end{align}
\end{proof}

\begin{remark}
Let $A(\cdot): [0,\tau] \to \mathcal{L}(D,X)$ be a Lipschitz continuous function, i.e.,
\begin{align*}
\|A(t)x-A(s)x\|\leq \lvert t-s \rvert \|x\|_D\qquad \left( t,s\in [0,\tau], x \in D \right).
\end{align*}
Assume that $A(t)\in\mathscr{MR}$ for any $t \in [0,\tau]$ and let $C \in \mathcal{L}(D,Y)$. By \eqref{estimate}, we obtain that if $\left(C,A(t)\right)$ is $p$-admissible for all $t \in [0,\tau]$, then $\left( C, A\left(\cdot\right) \right)$ is $p$-admissible for every $p \in (1,\infty)$.
\end{remark}

\begin{thm}\label{Kharou-thm2}
Let $A(\cdot): [0,\tau] \to \mathcal{L}(D,X)$ satisfies the condition {\bf(H)} and let $C \in \mathcal{L}(D,Y)$. Then $\left( C, A\left(\cdot\right) \right)$ is admissible, if and only if $\left(C,A(t)\right)$ is admissible for all $t \in [0,\tau)$.
\end{thm}

\begin{proof}
For $x \in D$ and $s \in [0,\tau)$, we have:
\begin{align*}
\int_{s}^{\tau} \| C U(t,s) x \|_{Y}^{p}  \mathrm{d}t &= \int_{s}^{\tau} \| C U(t,s) x - C e^{-(t-s)A(s)} x + C e^{-(t-s)A(s)} x \|_{Y}^{p} \mathrm{d}t \\
 &\leq 2^{p} \int_{s}^{\tau} \| C U(t,s) x - C e^{-(t-s)A(s)} x \|_{Y}^{p}  \mathrm{d}t \\
 &\qquad \qquad \qquad \qquad \qquad \qquad + 2^{p} \int_{s}^{\tau} \| C e^{-(t-s)A(s)} x \|_{Y}^{p} \mathrm{d}t \\
 &\leq 2^{p} \| C \|_{\mathcal{L}(D,Y)}^{p} \int_{s}^{\tau} \| U(t,s) x - e^{-(t-s)A(s)} x \|_{D}^{p}  \mathrm{d}t \\
 &\qquad \qquad \qquad \qquad \qquad \qquad + 2^{p} \int_{0}^{\tau-s} \| C e^{-t A(s)} x \|_{Y}^{p} \mathrm{d}t \\
 &\leq 2^{p} \| C \|_{\mathcal{L}(D,Y)}^{p} \| U(\cdot,s) x - e^{-(\cdot-s)A(s)} x \|_{L^{p}(s,\tau;D)}^{p} \\
 &\qquad \qquad \qquad \qquad \qquad \qquad + 2^{p} \gamma_{\tau -s}(s)^{p} \| x \|^{p} \\
 &\leq 2^{p} \| C \|_{\mathcal{L}(D,Y)}^{p} \| U(\cdot,a) x - e^{-(\cdot-s)A(s)} x \|_{\mathrm{MR}_{p}(s,\tau)}^{p} \\
 &\qquad \qquad \qquad \qquad \qquad \qquad + 2^{p} \gamma_{\tau-s}(s)^{p} \| x \|^{p} \\
 &\leq 2^{p} \| C \|_{\mathcal{L}(D,Y)}^{p} M^{p} \| x \|^{p} + 2^{p} \gamma_{\tau-s}(s)^{p} \| x \|^{p} ,
\end{align*}
where the last inequality follows from \eqref{key inequality}.

By the same arguments, we obtain the other implication. Let $x \in D$ and $[a,b]$ a sub-interval of $[0,\tau)$ with $a \neq b$, we have:
\begin{align*}
\int_{0}^{b-a} \left\| C e^{-t A(a)} x \right\|_{Y}^{p}  \mathrm{d}t &= \int_{a}^{b} \left\| C e^{-(t-a) A(a)} x \right\|_{Y}^{p}  \mathrm{d}t \\
 &= \int_{a}^{b} \left\| C e^{-(t-a) A(a)} x - C U(t,a)x + C U(t,a)x  \right\|_{Y}^{p}  \mathrm{d}t \\
 &\leq 2^p \int_{a}^{b} \left\| C U(t,a)x - C e^{-(t-a) A(a)} x \right\|_{Y}^{p}  \mathrm{d}t \\
 &\qquad \qquad \qquad \qquad \qquad \qquad+ 2^p \int_{a}^{b} \left\| C U(t,a)x \right\|_{Y}^{p}  \mathrm{d}t \\
 &\leq 2^{p} \left\| C \right\|_{\mathcal{L}(D,Y)}^{p} \left\| U(\cdot,a) x - e^{-(\cdot-a)A(a)} x \right\|_{\mathrm{MR}_{p}(a,b)}^{p} \\
 &\qquad \qquad \qquad \qquad \qquad \qquad + 2^p \gamma^p \left\| x \right\|^{p} \\
 &\leq 2^{p} \left\| C \right\|_{\mathcal{L}(D,Y)}^{p} M^p \| x \|^{p} + 2^p \gamma^p \| x \|^{p} .
\end{align*}
\end{proof}

\begin{remark}\label{H'}
If instead of {\bf(H)}, we make the following assumption {\bf(H')} on $A(\cdot)$:
\begin{itemize}
  \item [{\bf(H')}] For any $t\in [0,\tau],$ $A(t)\in\mathscr{MR}$, and 
  $$ \left\| (A(t)-A(s))e^{-(t-s)A(s)}x \right\| \leq \frac{c}{\lvert t-s \rvert^{\alpha}} \|x\| \qquad (t,s \in [0,\tau], x \in X) $$ for some constants $c > 0$ and $\alpha \in (0,1)$.
\end{itemize}
Then we can associate an evolution family on $X$ to $A(\cdot)$ and the result of Theorem \ref{Kharou-thm2} remains valid for $p \in (1,1/\alpha)$. The proof follows from the estimate $$ \left\| (A(t)-A(s))e^{-(t-s)A(s)}x \right\|^p \leq \frac{c^p}{\lvert t-s \rvert^{\alpha p}} \|x\|^p \qquad (t,s \in [0,\tau], x \in X) $$ for $p \in (1,1/\alpha)$. For instance, let $\mathbb{A} \in \mathcal{L}(D,X)$ such that $\mathbb{A} \in \mathscr{MR}$ and that the semigroup generated by $\mathbb{A}$ is of negative type. Consider the function $A(\cdot):[0,\tau] \to \mathcal{L}(D,X)$ given by $$ A(t)x = \mathbb{A} x + b(t) (- \mathbb{A})^{\alpha} x \qquad ( t \in [0,\tau], x \in D )$$ for some $\alpha \in (0,1)$, where $b:[0,\tau] \to [0,\infty)$ is a bounded and measurable map. Using Young's inequality and \cite[Example 2.9.(c)]{ACFP07}, it is easy to see that $A(t) \in \mathscr{MR}$ for any $t \in [0,\tau]$. Hence $A(\cdot)$ satisfies the condition {\bf(H')}.

\end{remark}

\section{Example}\label{sec:4}

Let $\Omega$ be a bounded domain in $\mathbb{R}^n$ with $n \geq 2$ with such that $\partial \Omega = \Gamma_0\cup \Gamma_1$ is bounded and of class $C^2$. We select 
\begin{align*}
 & X:=L^{2}(\Omega),\quad Y:=L^2(\partial \Omega), \cr & D:=\left\{ u\in H^{2}(\Omega) : u_{\lvert \Gamma_0}=\left( \frac{\partial u}{\partial \nu} \right)_{\lvert \Gamma_1}=0 \right\},
\end{align*}
where $\nu$ is the unit normal vector. Define $\mathbb{A}:=\Delta$ with domain $D(\mathbb{A})=D$. It was shown in \cite{LTZ04} that $\mathbb{A}$ is a generator of a strongly continuous semigroup on $X$. For $t \in [0,\tau]$, consider the operator $A(t):D\to X$ given by:
\begin{align*}
A(t)u= \mathbb{A} u+b(t) (- \mathbb{A})^{\alpha} u \qquad (u\in D),
\end{align*}
where $\alpha < 1/2$ and $b:[0,\tau] \to [0,\infty)$ is a continuous function. The operator $A(t)$ has maximal regularity for any $t \in [0,\tau]$. Hence $A(\cdot)$ satisfies the condition {\bf(H')} as in the remark above.\\
Consider the observation operator
\begin{align*}
C: D\to Y,\quad Cu=\begin{cases} 0,& \text{on}\,\Gamma_0,\cr u,& \text{on}\,\Gamma_1.\end{cases}
\end{align*}
According to \cite{LTZ04}, the operator $C$ is admissible for $\mathbb{A}$ and by \cite{HI06} $C$ is admissible for $A(t)$ for any $t \in [0,\tau]$. Finally, by Remark \ref{H'} (Theorem \ref{Kharou-thm2}), the operator $C$ is admissible for $A(\cdot)$.

\section*{Acknowledgments} 
I would like to thank Professor O. El-Mennaoui whose detailed comments helped me improve the organization and the content of the article.

\end{document}